

\baselineskip=14pt
\parskip=10pt
\def\halmos{\hbox{\vrule height0.15cm width0.01cm\vbox{\hrule height
  0.01cm width0.2cm \vskip0.15cm \hrule height 0.01cm width0.2cm}\vrule
  height0.15cm width 0.01cm}}

\magnification=\magstephalf

\def\1{{\overline{1}}}
\def\2{{\overline{2}}}
\parindent=0pt
\overfullrule=0in

\def\frac#1#2{{#1 \over #2}}


\centerline
{
\bf 
A Short Proof of a Ptolemy-Like Relation for an Even number of Points
}
\centerline
{\bf on a Circle Discovered by Jane McDougall}
\bigskip
\centerline
{
\it Marc CHAMBERLAND and Doron ZEILBERGER
}

\quad \quad \quad \quad \quad\quad \quad \quad
{\it In fond memory of Andrei Zelevinsky (1953-2013) who {\it loved} Ptolemy's theorem}

Jane McDougall[M] has discovered (and apparently found a complicated proof using heavy machinery) the
following beautiful Ptolemy-style theorem relating the distances 
amongst any even number of points on
a circle.

{\bf McDougall's Theorem}: Let $n$ be a positive integer, and let $P_i$ ($1 \leq i \leq 2n$) be
points on a circle. Let $d_{i,j}$ be the distance between $P_i$ and $P_j$, and let
$$
R_i:=\prod_{{ {1 \leq j \leq 2n} \atop {j \neq i} }}  d_{i,j} \quad .
$$
Then
$$
\sum_{i=1}^{n}  \frac{1}{R_{2i}} \, = \, \sum_{i=1}^{n}  \frac{1}{R_{2i-1}}  \quad .
$$

{\bf Proof}: Without loss of generality the circle is the unit circle. Let
$$
P_i=(\cos 2t_i \, , \, \sin 2t_i) \quad (1 \leq i \leq 2n) \quad.
$$
Thanks to trig, $d_{i,j}=2\sin  (t_j-t_i)$, and thanks to DeMoivre, this equals
$$
-\sqrt{-1} \left ( e^{\sqrt{-1} (t_j - t_i)}- e^{-\sqrt{-1} (t_j - t_i)} \right )  \quad.
$$
Let
$$
u_i=e^{\sqrt{-1} t_i} \quad .
$$
It follows that
$$
d_{i,j}= -\sqrt{-1} \left ( \frac{u_j}{u_i} - \frac{u_i}{u_j} \right ) = 
-\sqrt{-1} \frac{u_j^2-u_i^2}{u_i u_j} \quad (i<j) \quad .
$$
Hence
$$
(-\sqrt{-1})^{2n-1} 
(-1)^{i-1}\frac{1}{R_i}=(\prod_{j=1}^{2n} u_j) \frac{u_{i}^{2n-2}}{\prod_{j \neq i} (u_j^2-u_i^2) } \quad.
$$
So we have 
$$
(\sqrt{-1})^{2n-1} 
\sum_{i=1}^{2n}  \frac{(-1)^i}{R_{i}} \, = 
\left ( \prod_{j=1}^{2n} u_j\right ) 
\sum_{i=1}^{2n} \frac{u_{i}^{2n-2}}{\prod_{j \neq i} (u_j^2-u_i^2) } \quad.
\eqno(Jane)
$$

Recall the {\it Lagrange Interpolation Formula}: If $P(z)$ is a polynomial of degree $\leq N-1$ in $z$ then,
for any distinct numbers $z_1, \dots, z_N$,
$$
P(z)=\sum_{i=1}^{N}  
\frac{(z-z_1) \cdots (z-z_{i-1}) (z-z_{i+1}) \cdots (z-z_N)}
{(z_i-z_1) \cdots (z_i-z_{i-1}) (z_i-z_{i+1}) \cdots (z_i-z_N)} P(z_i) \quad .
\eqno(Joseph)
$$
(Let's recall the trivial proof: both sides are polynomials of degree $\leq N-1$ that coincide at the $N$ values
$z=z_1, \dots , z=z_N$, so they must be identically equal.)

Taking the polynomial to be $P(z)=z^{r}$ (with $r<N-1$), and 
equating the coefficient of $z^{N-1}$ on both sides of
$(Joseph)$, gives the identity:
$$
0=\sum_{i=1}^{N}  
\frac{z_i^{r}}
{(z_i-z_1) \cdots (z_i-z_{i-1}) (z_i-z_{i+1}) \cdots (z_i-z_N)}  \quad .
\eqno(Joseph')
$$
Now take $N=2n$, $r=n-1$, and $z_i=u_i^2$ in $(Joseph')$ and conclude that the right side of $(Jane)$ is indeed zero. \halmos 

It follows that if our $2n$ points lie on a line (the case $R=\infty$), the theorem is
true as well. Furthermore, in that case it is {\it even} true for an {\bf odd} number of points.
We leave this as an exercise to the dear readers.

{\bf Reference}

[M] Jane McDougall, {\it private communication} .

\bigskip
\hrule

\bigskip

Marc Chamberland, Department of Mathematics and Statistics, Grinnell College, Grinnell, Iowa 50112, USA.
{\tt chamberland at grinnell dot edu}

Doron Zeilberger, Mathematics Department, Rutgers University (New Brunswick), Piscataway, NJ 08854, USA.
{\tt zeilberg at math dot rutgers dot edu}

{\bf April 15, 2013}

\end